\documentclass[12pt,a4paper]{article}
\usepackage[english]{babel}
\usepackage{amsmath}
\usepackage{amssymb}
\usepackage{amsfonts}
\usepackage{hyperref}

\begin{document}


\renewcommand{\refname}{References}
\renewcommand{\contentsname}{Contents}

\begin{center}
{\huge Global Unique Solvability of the Initial-Boundary Value Problem for One-Dimensional Barotropic Equations of Viscous Compressible Bifluids}
\end{center}

\medskip

\begin{center}
{\large Alexander Mamontov,\quad Dmitriy Prokudin}
\end{center}

\medskip

\begin{center}
{\large March 20, 2020}
\end{center}

\medskip

\begin{center}
{
Lavrentyev Institute of Hydrodynamics, \\ Siberian Branch of the Russian Academy of Sciences\\
pr. Lavrent'eva 15, Novosibirsk 630090, Russia}
\end{center}

\medskip

\begin{center}
{\bfseries Abstract}
\end{center}


\begin{center}
\begin{minipage}{110mm}
We consider the equations of a multi-velocity model of a binary mixture of viscous compressible fluids (two-fluid medium) in the case of one-dimensional barotropic motions. We prove the global (in time) existence and uniqueness of a strong solution to the initial-boundary value problem describing the motion in a bounded spatial domain.
\end{minipage}
\end{center}

\bigskip


\newpage

\tableofcontents

\newpage

\section{Introduction}
\label{intro}

The paper is devoted to the problem of unique solvability for the equations of motion of viscous compressible fluid mixtures (multifluid media, multicomponent media, mul\-ti\-fuids). Regarding the origin of the model and its physical meaning, we refer to \cite{mamprok.muller}, \cite{mamprok.atkin}, \cite{mamprok.raj}, \cite{mamprok.nigm}. A survey of possible formulations of the model and known results can be found in \cite{mamprok.france}, \cite{mamprok.semi17}. Related multi-velocity models are discussed in \cite{mamprok.dorovski}, \cite{mamprok.dorov92}, \cite{mamprok.blokhin}, \cite{mamprok.gavr}, \cite{mamprok.jfcs2}. As the first results on the solvability of multifluid models in the multidimensional case (but in approximate formulations), we can indicate \cite{mamprok.frehsegm1}, \cite{mamprok.frehsegm2}, \cite{mamprok.frehsew}.

Weak solutions in the multidimensional case for the considered version of the model were constructed in \cite{mamprok.kucher11}, \cite{mamprok.kmp12}, \cite{mamprok.izvran14}. Similar results for related models were obtained in \cite{mamprok.semi162}, \cite{mamprok.jaimprok}, \cite{mamprok.smz1}, \cite{mamprok.semi161}, \cite{mamprok.smz2}, \cite{mamprok.izvran18}. For a number of reasons, including the goal of constructing more regular solutions, one-dimensional formulations are of interest. A detailed discussion of these formulations and an overview of the results can be found in \cite{mamprok.jfcs1}. In short, in multi-fluid models there are two dimensions, namely, the number of spatial variables and the number of components of the mixture. The equality of one of these parameters to one does not ``weaken'' the multidimensionality of the other parameter. Therefore, the study of one-dimensional flows of a multicomponent medium is of interest. Solvability for related one-dimensional models is shown in \cite{mamprok.semiprok17}, \cite{mamprok.jfcsprok}, \cite{mamprok.jms}.

The peculiarity of the paper is that a variant of the model with non-diagonal (namely, triangular) viscosity matrix is considered. Mathematically, we are talking about a system of equations of a mixed type, obtained by combining compressible Navier---Stokes-type systems, but additionally interconnected in the higher order terms due to the viscosity matrix and therefore not allowing the automatic transfer of the known results of the classical (single-fluid) theory developed in such works as \cite{mamprok.kazhshel}, \cite{mamprok.akm}, \cite{mamprok.shel84}, \cite{mamprok.kazhsel}.

Solvability results for a model with a diagonal viscosity matrix are obtained in \cite{mamprok.kazhpetr}, \cite{mamprok.petr82}. However, it should be noted that we do not use here the technique of these works, but we directly extend the technique developed for a single-component viscous gas to the multi-component case with a non-diagonal viscosity matrix.

\section{Statement of the Problem and the Main Result}
\label{sec1}

\noindent\indent We consider the system of equations of motion of binary viscous compressible fluid mixtures in the case of one spatial variable $(i=1,2)$:
\begin{equation}\label{postzad_1}
\partial_{t}\rho_{i}+\partial_{x}(\rho_{i} u_{i})=0,
\end{equation}
\begin{equation}\label{postzad_2}
\rho_{i}\left(\partial_{t}u_{i}+u_{i}\partial_{x}u_{i}\right)+K_{i}\partial_{x}(\rho^{\gamma_{i}}_{i})=
\sum\limits_{j=1}^2 \mu_{ij}\partial_{xx}u_{j}+(-1)^{i+1}a(u_{2}-u_{1}).
\end{equation}
Here $\rho_{i}$ (the density of the $i$-th constituent of the mixture) and $u_{i}$ (the velocity of the $i$-th component) are unknown values, and we are given physical constants $a>0$, $K_{i}>0$,  $\gamma_{i}>1$, $i=1,2$, and constant viscosity coefficients $\{\mu_{ij}\}_{i, j=1}^{2}$, which form the matrix $\textbf{M}>0$, i.~e.
$(\textbf{M}\boldsymbol{\xi},\boldsymbol{\xi})\geqslant M_{0}|\boldsymbol{\xi}|^2$ for all $\boldsymbol{\xi}\in{\mathbb R}^{2}$ with some constant $M_{0}(\textbf{M})>0$.

It should be noted that the authors do not claim the maximum generality of the equation of state and structure of momentum exchange terms (which was achieved, for example, in
\cite{mamprok.zlot95}), since the main purpose of the article is to focus on the analysis of the viscous terms.

Let us consider the system \eqref{postzad_1}, \eqref{postzad_2} in the rectangular $Q_{T}$ (here and below $Q_{t}=(0, 1)\times(0, t)$) with an arbitrary finite positive height $T$ complemented by the following initial and boundary conditions ($i=1,\, 2$):
\begin{equation}\label{postzad_3}
\rho_{i}|_{t=0}=\rho_{0i}(x), \quad u_{i}|_{t=0}=u_{0i}(x),\quad x\in [0, 1],
\end{equation}
\begin{equation}\label{postzad_4}
u_{i}|_{x=0}=u_{i}|_{x=1}=0,\quad t\in [0, T].
\end{equation}

\medskip

{\bfseries Definition 1.} {\it By a strong solution to the problem \eqref{postzad_1}--\eqref{postzad_4} we mean a collection of functions $(\rho_{1}, \rho_{2}, u_{1}, u_{2})$ such that the equations \eqref{postzad_1}, \eqref{postzad_2} are valid almost everywhere in $Q_{T}$, the initial data \eqref{postzad_3} are accepted for almost all $x\in (0, 1)$, the boundary conditions \eqref{postzad_4} are valid for a.~a. $ t\in (0, T)$, and the following inequalities and inclusions hold $(i=1,\, 2)$
\begin{equation}\label{postzad_5}
\begin{array}{c}
\displaystyle
\rho_{i}\geqslant {\rm const}>0\text{ a. e. in }Q_T,
\quad \rho_{i}\in L_{\infty}\big(0, T; W^{1}_{2}(0, 1)\big),\\\\  \partial_{t}\rho_{i}\in L_{\infty}\big(0, T; L_{2}(0, 1)\big),\\\\
\displaystyle u_{i}\in L_{\infty}\big(0, T; W^{1}_{2}(0, 1)\big)\bigcap L_{2}\big(0, T; W^{2}_{2}(0, 1)\big),\\\\
\displaystyle \partial_{t}u_{i} \in L_{2}(Q_{T});
\end{array}
\end{equation}
thus, \eqref{postzad_1} are understood as equalities in $L_{\infty}\big(0, T; L_{2}(0, 1)\big)$, and \eqref{postzad_2}\linebreak in~$L_{2}(Q_{T})$.}

\medskip

{\bfseries Theorem 1.} {\it Let the initial data in \eqref{postzad_3} satisfy the conditions
\begin{equation}\label{postzad_6}
\rho_{0i}\in W^{1}_{2}(0,1),\quad \rho_{0i}>0,\quad u_{0i}\in \overset{\circ}{W^1_2}(0, 1),\quad i=1,\, 2,
\end{equation}
and $\mu_{12}=0$. Then there exists the unique strong solution to the problem \eqref{postzad_1}---\eqref{postzad_4} in the sense of Definition~1.}

\medskip

{\bfseries Proof.} The existence of a unique strong solution to the problem~\eqref{postzad_1}---\eqref{postzad_4} in a small time interval $[0, t_{0}]$ is proved using the Galerkin method, see \cite{mamprok.jms}. To extend the solution from the segment $[0, t_{0}]$ to the entire segment $[0, T]$ under consideration, it is necessary to obtain a priori estimates for the local solution in the classes \eqref{postzad_5}, the limiting constants in which depend only on the input data and on the value of $T$, but not on the parameter $t_{0}$. Therefore, we will focus on obtaining global estimates.

Consider a hypothetical solution $(\rho_{1}, \rho_{2}, u_{1}, u_{2}) $ of the problem \eqref{postzad_1}---\eqref{postzad_4}, which has all the necessary differential properties, and such that the densities $\rho_{i}$, $i = 1,2$ are strictly positive and bounded (see Definition 1).

\section{Lagrangian Mass Coordinates}
\label{sec2}

\noindent\indent We transform the problem \eqref{postzad_1}--\eqref{postzad_4} introducing two variants of the Lagrangian coordinates associated with the velocity of each component of the mixture, with the number $m\in\{1, 2\}$. In other words, we take\linebreak
$\displaystyle y_{m}(x,t)=\int\limits_{0}^{x}\rho_{m}(s,t)\,ds$ and $t$ as new independent variables. Then the system \eqref{postzad_1}, \eqref{postzad_2} takes the form
\begin{equation}\label{postzad_7}
\partial_{t}\rho_{i}+\rho_{m}(u_{i}-u_{m})\partial_{y_{m}}\rho_{i}+\rho_{i}\rho_{m}\partial_{y_{m}}u_{i}=0,\quad i,m=1,2,
\end{equation}
\begin{equation}\label{postzad_8}
\begin{array}{c}
\displaystyle
\frac{\rho_{i}}{\rho_{m}}\partial_{t}u_{i}+\rho_{i}(u_{i}-u_{m})\partial_{y_{m}}u_{i}+K_{i}\partial_{y_{m}} \rho_{i}^{\gamma_{i}}=
\sum\limits_{j=1}^2 \mu_{ij}\partial_{y_{m}}(\rho_{m}\partial_{y_{m}}u_{j})+ \\ \\
\displaystyle
+\frac{(-1)^{i+1}a}{\rho_{m}}(u_{2}-u_{1}),\quad i, m=1, 2.
\end{array}
\end{equation}
The domain $Q_{T}$, during this transformation, maps into the rectangular\linebreak $(0, d_{m})\times(0, T)$, where $\displaystyle d_{m}=\int\limits_{0}^{1}\rho_{0m}\,dx>0$, initial and boundary conditions are converted to the form (here $i,m=1, 2$)
\begin{equation}\label{postzad_9}
\rho_{i}|_{t=0}=\widetilde{\rho}_{0i}(y_{1}), \quad u_{i}|_{t=0}=\widetilde{u}_{0i}(y_{1}),\quad y_{1}\in [0, d_{1}]
\end{equation}
(in the first Lagrangian coordinate system), or
\begin{equation}\label{postzad_10}
\rho_{i}|_{t=0}=\widetilde{\widetilde{\rho}}_{0i}(y_{2}), \quad u_{i}|_{t=0}=\widetilde{\widetilde{u}}_{0i}(y_{2}),\quad y_{2}\in [0, d_{2}]
\end{equation}
(in the second one), and in both cases
\begin{equation}\label{postzad_11}
u_{i}|_{y_{m}=0}=u_{i}|_{y_{m}=d_{m}}=0,\quad t\in [0, T].
\end{equation}
Let us note that the equations \eqref{postzad_7} and \eqref{postzad_8} lead to the relations
$$\partial_{t}\left(\frac{\rho_{i}}{\rho_{m}}\right)+\partial_{y_{m}}(\rho_{i}(u_{i}-u_{m}))=0,\quad i, m= 1,2,$$
$$\begin{array}{c}\displaystyle
\partial_{t}\left(\frac{\rho_{i}}{\rho_{m}}u_{i}\right)+\partial_{y_{m}}(\rho_{i}(u_{i}-u_{m})u_{i})+K_{i}\partial_{y_{m}} \rho_{i}^{\gamma_{i}}=
\sum\limits_{j=1}^2 \mu_{ij}\partial_{y_{m}}(\rho_{m}\partial_{y_{m}}u_{j})+ \\ \\
\displaystyle
+\frac{(-1)^{i+1}a}{\rho_{m}}(u_{2}-u_{1})\quad i, m= 1,2.
\end{array}$$

\medskip

{\bfseries Remark 1.} {\it After introducing the Lagrangian coordinates associated with the velocity of each component, the original $($one$)$ system of equations \eqref{postzad_1}, \eqref{postzad_2} is transformed into two systems \eqref{postzad_7}, \eqref {postzad_8}, and each of the systems is defined in its own spatial region, which significantly complicates their analysis and, therefore, a simplifying assumption was made on the triangularity of the viscosity matrix. However this statement of the problem is an essential advance in the theory, since $($as it was mentioned above$)$ so far, the results were obtained only for the diagonal matrix of viscosities \cite{mamprok.kazhpetr}, \cite{mamprok.petr82}.}

\medskip

{\bfseries Remark 2.} {\it The systems \eqref{postzad_7} and \eqref{postzad_8} are of divergent form as $i\neq m$.}

\medskip

Further, in Sections 4--6, we reproduce a standard technique for estimating the so\-lu\-tions of one-dimensional viscous gas equations with small variations according to bifluid specificity. However, this work needs to be done. Starting from Section~7, the differences from the one-component case become more significant.

\section{First A Priori Estimate}
\label{sec3}

\noindent\indent Let us multiply the equations \eqref{postzad_2} by $ u_{i}$, integrate over $(0, 1)$ and sum in $i$ from $1$ to $2$. Taking into account that, due to \eqref{postzad_1}, \eqref{postzad_4} and the condition $\textbf{M}>0$, the following relations hold
\begin{equation}\label{step1_1}
\sum\limits_{i=1}^{2}\int\limits^{1}_{0}\Big(\rho_{i}\partial_{t}u_{i}+\rho_{i} u_{i}\partial_{x}u_{i}\Big)u_{i}\,
dx=\frac{d}{dt}\left(\sum\limits_{i=1}^{2}\int\limits^{1}_{0}\frac{1}{2}\rho_{i} u_{i}^{2}\, dx\right),
\end{equation}
\begin{equation}\label{step1_2}
\begin{array}{c}
\displaystyle
\sum\limits_{i=1}^{2}\int\limits^{1}_{0}u_{i}K_{i}\partial_{x}\rho_{i}^{\gamma_{i}}\,
dx=-\sum\limits_{i=1}^{2}\int\limits^{1}_{0}K_{i}\rho_{i}^{\gamma_{i}}(\partial_{x}u_{i})\, dx=\\\\
\displaystyle
=\frac{d}{dt}\left(\sum\limits_{i=1}^{2}\int\limits^{1}_{0}\frac{K_{i}}{\gamma_{i}-1}\rho_{i}^{\gamma_{i}}\, dx\right),
\end{array}
\end{equation}
\begin{equation}\label{step1_3}
\begin{array}{c}
\displaystyle
\sum\limits_{i,j=1}^2
\mu_{ij}\int\limits\limits_{0}^{1}(\partial_{xx}u_{j})u_{i}\,dx=-\sum\limits_{i,j=1}^2
\mu_{ij}\int\limits^{1}_{0}(\partial_{x}u_{i})(\partial_{x}u_{j})\,dx\leqslant\\\\
\displaystyle
\leqslant -M_{0}
\sum\limits_{i=1}^2\int\limits\limits_{0}^{1}|\partial_{x}u_{i}|^{2}\,dx,
\end{array}
\end{equation}
we arrive at the inequality
\begin{equation}\label{step1_4}
\begin{array}{c}
\displaystyle
\frac{d}{dt}\left(\sum\limits_{i=1}^{2}
\int\limits^{1}_{0}\left(\frac{1}{2}\rho_{i} u_{i}^{2}+\frac{K_{i}}{\gamma_{i}-1}\rho_{i}^{\gamma_{i}}\right)\,dx\right)+
M_{0}\sum\limits_{i=1}^2\int\limits\limits_{0}^{1}|\partial_{x}u_{i}|^{2}\,dx+\\\\
\displaystyle +a\int\limits^{1}_{0}|u_{1}-u_{2}|^{2}\,dx\leqslant 0.
\end{array}
\end{equation}
Integrating the inequality \eqref{step1_4} over $(0, t)$, we obtain, using \eqref{postzad_3}, the first a priori estimate
\begin{equation}\label{step1_5}\begin{array}{c}\displaystyle
\sum\limits_{i=1}^{2}\int\limits^{1}_{0}\left(\frac{1}{2}\rho_{i} u_{i}^{2}+\frac{K_{i}}{\gamma_{i}-1}\rho_{i}^{\gamma_{i}}\right)\,dx+
M_{0}\sum\limits_{i=1}^2\int\limits\limits_{0}^{t}\int\limits\limits_{0}^{1}|\partial_{x}u_{i}|^{2}\,dxd\tau+\\ \\
\displaystyle
+a\int\limits\limits_{0}^{t}\int\limits^{1}_{0}|u_{1}-u_{2}|^{2}\,dxd\tau\leqslant\sum\limits_{i=1}^{2}
\int\limits^{1}_{0}\left(\frac{1}{2}\rho_{0i}u_{0i}^{2}+\frac{K_{i}}{\gamma_{i}-1}\rho_{0i}^{\gamma_{i}}\right)\,dx=:B_{1}.
\end{array}
\end{equation}
The estimate \eqref{step1_5} obviously provide the inequalities
\begin{equation}\label{step2_17}
\left\|u_{i}\right\|_{L_{2}\big(0, T;L_{\infty}(0, 1)\big)}\leqslant\sqrt{\frac{B_{1}}{M_{0}}},\quad i=1, 2.
\end{equation}
In the Lagrangian coordinates $(y_{m}, t)$ the formula \eqref{step1_5} takes the form
\begin{equation}\label{step1_6}\begin{array}{c}\displaystyle
\sum\limits_{i=1}^{2}
\int\limits^{d_{m}}_{0}\left(\frac{1}{2}\frac{\rho_{i}u_{i}^{2}}{\rho_{m}}+\frac{K_{i}}{\gamma_{i}-1}\frac{\rho_{i}^{\gamma_{i}}}{\rho_{m}}\right)\,dy_{m}+
M_{0}\sum\limits_{i=1}^2\int\limits\limits_{0}^{t}\int\limits\limits_{0}^{d_{m}}\rho_{m}|\partial_{y_{m}}u_{i}|^{2}\,dy_{m}d\tau+\\ \\
\displaystyle
+a\int\limits\limits_{0}^{t}\int\limits^{d_{m}}_{0}\frac{|u_{1}-u_{2}|^{2}}{\rho_{m}}\,dy_{m}d\tau\leqslant B_{1},\quad t\in [0, T],\quad m= 1,2,
\end{array}
\end{equation}
from which, in particular, we have
\begin{equation}\label{step1_7}\begin{array}{c}\displaystyle
\int\limits^{d_{m}}_{0}\left(\frac{u_{m}^{2}}{2}+\frac{K_{m}}{\gamma_{m}-1}\rho_{m}^{\gamma_{m}-1}\right)\,dy_{m}+
M_{0}\int\limits\limits_{0}^{t}\int\limits\limits_{0}^{d_{m}}\rho_{m}|\partial_{y_{m}}u_{m}|^{2}\,dy_{m}d\tau+\\ \\
\displaystyle
+a\int\limits\limits_{0}^{t}\int\limits^{d_{m}}_{0}\frac{|u_{1}-u_{2}|^{2}}{\rho_{m}}\,dy_{m}d\tau\leqslant B_{1},\quad t\in [0, T],\quad m= 1,2.
\end{array}
\end{equation}

\section{Estimate of Strict Positivity and\\ Boundedness of the Density $\rho_{1}$}
\label{sec4}

\noindent\indent In the system \eqref{postzad_7}, \eqref{postzad_8}, we consider the equations with the number\linebreak $i=m=1$:
\begin{equation}\label{step2_1}
\partial_{t}\rho_{1}+\rho_{1}^{2}\partial_{y_{1}}u_{1}=0,
\end{equation}
\begin{equation}\label{step2_2}
\partial_{t}u_{1}+K_{1}\partial_{y_{1}} \rho_{1}^{\gamma_{1}}=\mu_{11}\partial_{y_{1}}(\rho_{1}\partial_{y_{1}}u_{1})+\frac{a}{\rho_{1}}(u_{2}-u_{1}).
\end{equation}
From the equation \eqref{step2_1} we express the value
\begin{equation}\label{step2_3}
\rho_{1}\partial_{y_{1}}u_{1}=-\partial_{t}\ln{\rho_{1}}
\end{equation}
and substitute into \eqref{step2_2}:
\begin{equation}\label{step2_4}
\mu_{11}\partial_{ty_{1}}\ln{\rho_{1}}+K_{1}\partial_{y_{1}} \rho_{1}^{\gamma_{1}}=-\partial_{t}u_{1}+\frac{a}{\rho_{1}}(u_{2}-u_{1}).
\end{equation}
We multiply the equality \eqref{step2_4} by $\partial_{y_{1}}\ln{\rho_{1}}$ and integrate in $y_{1}$ from $0$ to $d_{1}$, then we get
\begin{equation}\label{step2_5}
\begin{array}{c}
\displaystyle
\frac{\mu_{11}}{2}\frac{d}{dt}\int\limits^{d_{1}}_{0}|\partial_{y_{1}}\ln{\rho_{1}}|^{2}\,
dy_{1}+K_{1}\gamma_{1}\int\limits^{d_{1}}_{0}\rho_{1}^{\gamma_{1}}|\partial_{y_{1}}\ln{\rho_{1}}|^{2}\,dy_{1}=\\ \\
\displaystyle
=-\int\limits^{d_{1}}_{0}(\partial_{t}u_{1})(\partial_{y_{1}}\ln{\rho_{1}})\,
dy_{1}+a\int\limits^{d_{1}}_{0}\frac{u_{2}-u_{1}}{\rho_{1}}(\partial_{y_{1}}\ln{\rho_{1}})\,dy_{1}.
\end{array}
\end{equation}
The first term in the right-hand side of \eqref{step2_5} is transformed by integrating by parts and taking into account \eqref{postzad_11}, \eqref{step2_3}:
\begin{equation}\label{step2_6}
-\int\limits^{d_{1}}_{0}(\partial_{t}u_{1})(\partial_{y_{1}}\ln{\rho_{1}})\,
dy_{1}=-\frac{d}{dt}\int\limits^{d_{1}}_{0}u_{1}(\partial_{y_{1}}\ln{\rho_{1}})\,
dy_{1}+\int\limits^{d_{1}}_{0}\rho_{1}|\partial_{y_{1}}u_{1}|^{2}\,dy_{1},
\end{equation}
and the second one can be estimated from above:
\begin{equation}\label{step2_7}
\begin{array}{c}
\displaystyle
a\int\limits^{d_{1}}_{0}\frac{u_{2}-u_{1}}{\rho_{1}}(\partial_{y_{1}}\ln{\rho_{1}})\,
dy_{1}\leqslant \\\\
\displaystyle
\leqslant
a\|\rho_{1}^{-1/2}\|_{C[0, d_{1}]}\|\partial_{y_{1}}\ln{\rho_{1}}\|_{L_{2}(0, d_{1})}\|\rho_{1}^{-1/2}(u_{1}-u_{2})\|_{L_{2}(0, d_{1})}.
\end{array}
\end{equation}
From \eqref{postzad_1}, \eqref{postzad_3} and \eqref{postzad_4} it follows clearly that for each $t\in[0,T]$, in at least one point  $\varsigma_{1}(t)\in[0, d_{1}]$, the relation
\begin{equation}\label{step2_8}
\rho_{1}(\varsigma_{1}(t), t)=d_{1}
\end{equation}
holds. Thus, the representation
\begin{equation}\label{step2_9}
\rho_{1}^{-1/2}(y_{1},t)=d_{1}^{-1/2}+\int\limits_{\varsigma_{1}}^{y_{1}}\partial_{s}\rho_{1}^{-1/2}(s,t)\, ds
\end{equation}
is valid, from which by the Cauchy inequality, we have
\begin{equation}\label{step2_10}
\|\rho_{1}^{-1/2}\|_{C[0, d_{1}]}\leqslant d_{1}^{-1/2}+\frac{1}{2}\|\partial_{y_{1}}\ln{\rho_{1}}\|_{L_{2}(0, d_{1})}.
\end{equation}
Applying the inequalities \eqref{step2_7} and \eqref{step2_10}, after integration of \eqref{step2_5} in $t$ we get
\begin{equation}\label{step2_11}
\begin{array}{c}
\displaystyle
\mu_{11}\|\partial_{y_{1}}\ln{\rho_{1}}\|^{2}_{L_{2}(0, d_{1})}+2K_{1}\gamma_{1}\int\limits_{0}^{t}\int\limits^{d_{1}}_{0}\rho_{1}^{\gamma_{1}}|\partial_{y_{1}}\ln{\rho_{1}}|^{2}\,
dy_{1}d\tau\leqslant\\ \\
\displaystyle
\leqslant \mu_{11}\|\partial_{y_{1}}\ln{\widetilde{\rho}_{01}}\|^{2}_{L_{2}(0, d_{1})} -2\int\limits^{d_{1}}_{0}u_{1}(\partial_{y_{1}}\ln{\rho_{1}})\,
dy_{1}+\\ \\
\displaystyle
+2\int\limits^{d_{1}}_{0}\widetilde{u}_{01}(\partial_{y_{1}}\ln{\widetilde{\rho}_{01}})\,
dy_{1}+2\int\limits_{0}^{t}\int\limits^{d_{1}}_{0}\rho_{1}|\partial_{y_{1}}u_{1}|^{2}\,dy_{1}d\tau+\\ \\
\displaystyle
+a\int\limits_{0}^{t}\|\rho_{1}^{-1/2}(u_{1}-u_{2})\|_{L_{2}(0, d_{1})}\Big(2d_{1}^{-1/2}+\\\\
\displaystyle
+\|\partial_{y_{1}}\ln{\rho_{1}}\|_{L_{2}(0, d_{1})}\Big)\|\partial_{y_{1}}\ln{\rho_{1}}\|_{L_{2}(0, d_{1})}\,d\tau.
\end{array}
\end{equation}
Using the Cauchy inequality
$$
-2\int\limits^{d_{1}}_{0}u_{1}(\partial_{y_{1}}\ln{\rho_{1}})\,
dy_{1}\leqslant \frac{\mu_{11}}{2}\|\partial_{y_{1}}\ln{\rho_{1}}\|^{2}_{L_{2}(0, d_{1})}+\frac{2}{\mu_{11}}\|u_{1}\|^{2}_{L_{2}(0, d_{1})}
$$
and the estimate \eqref{step1_7}, from \eqref{step2_11} we deduce
\begin{equation}\label{step2_12}
\begin{array}{c}
\displaystyle
\|\partial_{y_{1}}\ln{\rho_{1}}\|^{2}_{L_{2}(0, d_{1})}\leqslant B_{2}+\\\\
\displaystyle
+B_{3}\int\limits_{0}^{t}\|\rho_{1}^{-1/2}(u_{1}-u_{2})\|_{L_{2}(0, d_{1})}\|\partial_{y_{1}}\ln{\rho_{1}}\|^{2}_{L_{2}(0, d_{1})}\,d\tau,
\end{array}
\end{equation}
where $B_{2}=B_{2}\left(M_{0}, B_{1}, \mu_{11}, a, T, d_{1}, \|\partial_{y_{1}}\ln{\widetilde{\rho}_{01}}\|_{L_{2}(0, d_{1})}, \|\widetilde{u}_{01}\|_{L_{2}(0, d_{1})}\right)$,\linebreak $\displaystyle B_{3}=\frac{4a}{\mu_{11}}$. Using \eqref{step1_7} again, by the Gronwall lemma from \eqref{step2_12} we obtain
\begin{equation}\label{step2_13}
\|\partial_{y_{1}}\ln{\rho_{1}}(t)\|_{L_{2}(0, d_{1})}\leqslant B_{4}\quad \forall\, t\in [0, T],
\end{equation}
where $B_{4}=B_{4}(B_{1}, B_{2}, B_{3}, a, T)$. Hence,
\begin{equation}\label{step2_14}\begin{array}{c}\displaystyle
|\ln\rho_{1}(y_{1},t)|=\left|\ln\rho_{1}(\varsigma_{1},t)+\int\limits_{\varsigma_{1}}^{y_{1}}\partial_{s}\ln\rho_{1}(s,t)\, ds\right|\leqslant \\ \\
\displaystyle
\leqslant |\ln{d_{1}}|+\sqrt{d_{1}}\|\partial_{y_{1}}\ln{\rho_{1}}\|_{L_{2}(0, d_{1})}\leqslant B_{5}(B_{4}, d_{1}),
\end{array}
\end{equation}
and consequently
\begin{equation}\label{step2_15}
0<B_{6}^{-1}\leqslant\rho_{1}(y_{1}, t)\leqslant B_{6}<\infty\quad\text{as}\quad(y_{1},t)\in[0, d_{1}]\times[0, T],
\end{equation}
where $B_{6}=e^{B_{5}}$. Finally,
\begin{equation}\label{step2_16}
0<B_{6}^{-1}\leqslant\rho_{1}(x, t)\leqslant B_{6}<\infty\quad\text{as}\quad(x,t)\in[0, 1]\times[0, T].
\end{equation}

\section{Estimates for the Derivatives of $\rho_{1}$ and $u_{1}$}
\label{sec5}

\noindent\indent After the boundedness and positivity of the density $\rho_{1} $ is shown, a priori estimates for the derivatives of $\rho_{1} $ and $u_{1} $ can be obtained in the Eulerian variables $(x, t)$. So, for example, from \eqref{step2_13} and \eqref{step2_16} we deduce
\begin{equation}\label{step3_1}
\|\partial_{x}\rho_{1}(t)\|_{L_{2}(0, 1)}\leqslant B_{7}=B_{4}B_{6}^{\frac{3}{2}}\quad \forall\, t\in [0, T].
\end{equation}
Next, we square the equation \eqref{postzad_2} with the number $i=1$, and obtain
\begin{equation}\label{step3_2}
\begin{array}{c}
\displaystyle
\rho_{1}|\partial_{t}u_{1}|^{2}+\frac{\mu_{11}^{2}}{\rho_{1}}|\partial_{xx}u_{1}|^{2}-2\mu_{11}(\partial_{t}u_{1})(\partial_{xx}u_{1})=\\\\
\displaystyle=\frac{1}{\rho_{1}}\left|K_{1}\partial_{x}\rho_{1}^{\gamma_{1}}+\rho_{1}u_{1}\partial_{x}u_{1}+a(u_{1}-u_{2})\right|^{2}.
\end{array}
\end{equation}
Let us introduce the function $\alpha(t)$ as
$$
\alpha(t)=\mu_{11}\int\limits_{0}^{1}|\partial_{x}u_{1}|^{2}\, dx+\int\limits_{0}^{t}\int\limits_{0}^{1}\left(\rho_{1}|\partial_{t}u_{1}|^{2}+\frac{\mu_{11}^2}{\rho_{1}}|\partial_{xx}u_{1}|^{2}\right)\, dxd\tau.
$$
From \eqref{postzad_4}, \eqref{step3_2} and the inequalities \eqref{step2_16}, \eqref{step3_1} it follows that
\begin{equation}\label{step3_3}
\alpha^{\prime}(t)\leqslant B_{8}+B_{9}\|u_{1}-u_{2}\|^{2}_{L_{2}(0,1)}+B_{10}\|u_{1}\|^{2}_{L_{\infty}(0,1)}\alpha(t),
\end{equation}
where $B_{8}=B_{8}(B_{6}, B_{7}, K_{1}, \gamma_{1})$, $B_{9}=3a^{2}B_{6}$, $\displaystyle B_{10}=\frac{3B_{6}}{\mu_{11}}$, and hence due to \eqref{step1_5}, \eqref{step2_17} by the Gronwall lemma the inequality
\begin{equation}\label{step3_4}
\alpha(t)\leqslant B_{11}\left(B_{1}, B_{8}, B_{9}, B_{10}, M_{0}, \mu_{11}, a, T, \|\partial_{x}u_{01}\|_{L_{2}(0, 1)}\right)
\end{equation}
follows, from which we arrive at the estimate
\begin{equation}\label{step3_5}
\begin{array}{c}
\displaystyle
\|\partial_{x}u_{1}\|_{L_{\infty}\left(0, T; L_{2}(0, 1)\right)}+\|\partial_{xx}u_{1}\|_{L_{2}(Q_{T})}+\\\\
\displaystyle
+\|\partial_{t}u_{1}\|_{L_{2}(Q_{T})}\leqslant B_{12}(B_{6}, B_{11}, \mu_{11}).
\end{array}
\end{equation}
Then from the continuity equation \eqref{postzad_1} with the number $i = 1 $ and the estimates \eqref{step2_16}, \eqref{step3_1}, \eqref{step3_5} it follows that
\begin{equation}\label{step3_6}
\|\partial_{t}\rho_{1}\|_{L_{\infty}(0, T; L_{2}(0,1))}\leqslant B_{13}(B_{6}, B_{7}, B_{12}).
\end{equation}

\section{Estimate of Strict Positivity and\\ Boundedness of the Density $\rho_{2}$}
\label{sec6}

\noindent\indent Let us consider the equations with the number $i=m=2$ in the system \eqref{postzad_7},  \eqref{postzad_8}:
\begin{equation}\label{step4_1}
\partial_{t}\rho_{2}+\rho_{2}^{2}\partial_{y_{2}}u_{2}=0,
\end{equation}
\begin{equation}\label{step4_2}
\partial_{t}u_{2}+K_{2}\partial_{y_{2}} \rho_{2}^{\gamma_{2}}=
 \mu_{21}\partial_{y_{2}}(\rho_{2}\partial_{y_{2}}u_{1})+\mu_{22}\partial_{y_{2}}(\rho_{2}\partial_{y_{2}}u_{2})
-\frac{a}{\rho_{2}}(u_{2}-u_{1}).
\end{equation}

\medskip

{\bfseries Remark 3.} {\it The equations \eqref{step2_2} and \eqref{step4_2} significantly differ because the first term in the right-hand side of \eqref{step4_2} has a nontrivial form due to the cross product $\rho_{2} \partial_{y_{2}} u_{1} $, and therefore the estimate for $\rho_{2} $ cannot be obtained in the same way as it was done for $ \rho_{1} $ $($in Section 5$)$. Two variants of the analysis of \eqref{step4_2} arise here: either to deal with this nontrivial term as a highest derivative of the sought function, or to consider it as a right-hand side term. What to do in the first case is not obvious yet, perhaps we should use the equations \eqref{postzad_7}, \eqref{postzad_8} with $ i \neq m $, but how exactly to do this is not clear $($if this difficulty is overcome, general viscosity matrix problem will be solved$)$. The second option is not easy either, because, on the one hand, we already have estimates for $u_{1} $, but on the other hand, $ \rho_{2} $, which is to be estimated, is present in the nontrivial term. We managed to overcome this difficulty $($as will be clear below$)$.}

\medskip

Let us express from \eqref{step4_1}
\begin{equation}\label{step4_3}
\rho_{2}\partial_{y_{2}}u_{2}=-\partial_{t}\ln{\rho_{2}}
\end{equation}
and substitute into \eqref{step4_2}:
\begin{equation}\label{step4_4}
\mu_{22}\partial_{ty_{2}}\ln{\rho_{2}}+K_{2}\partial_{y_{2}} \rho_{2}^{\gamma_{2}}=
 -\partial_{t}u_{2}
+\frac{a}{\rho_{2}}(u_{1}-u_{2})+\mu_{21}\partial_{y_{2}}(\rho_{2}\partial_{y_{2}}u_{1}).
\end{equation}
Let us multiply \eqref{step4_4} by $\partial_{y_{2}}\ln{\rho_{2}}$ and integrate in $y_{2}$ from $0$ to $d_{2}$, then we obtain
\begin{equation}\label{step4_5}\begin{array}{c}\displaystyle
\frac{\mu_{22}}{2}\frac{d}{dt}\int\limits^{d_{2}}_{0}|\partial_{y_{2}}\ln{\rho_{2}}|^{2}\,
dy_{2}+K_{2}\gamma_{2}\int\limits^{d_{2}}_{0}\rho_{2}^{\gamma_{2}}|\partial_{y_{2}}\ln{\rho_{2}}|^{2}\,
dy_{2}=\\\\
\displaystyle
=-\int\limits^{d_{2}}_{0}(\partial_{t}u_{2})(\partial_{y_{2}}\ln{\rho_{2}})\,dy_{2}+a\int\limits^{d_{2}}_{0}\frac{u_{1}-u_{2}}{\rho_{2}}(\partial_{y_{2}}\ln{\rho_{2}})\,
dy_{2}+\\ \\
\displaystyle
+\mu_{21}\int\limits^{d_{2}}_{0}(\partial_{y_{2}}(\rho_{2}\partial_{y_{2}}u_{1}))(\partial_{y_{2}}\ln{\rho_{2}})\,dy_{2}.
\end{array}
\end{equation}
The first term in the right-hand side of \eqref{step4_5} is transformed by integrating by parts and taking into account \eqref{postzad_11}, \eqref{step4_3}:
\begin{equation}\label{step4_6}
-\int\limits^{d_{2}}_{0}(\partial_{t}u_{2})(\partial_{y_{2}}\ln{\rho_{2}})\,
dy_{2}=-\frac{d}{dt}\int\limits^{d_{2}}_{0}u_{2}(\partial_{y_{2}}\ln{\rho_{2}})\,
dy_{2}+\int\limits^{d_{2}}_{0}\rho_{2}|\partial_{y_{2}}u_{2}|^{2}\,dy_{2},
\end{equation}
and all other terms are estimated from above:
\begin{equation}\label{step4_7}
\begin{array}{c}
\displaystyle
a\int\limits^{d_{2}}_{0}\frac{u_{1}-u_{2}}{\rho_{2}}(\partial_{y_{2}}\ln{\rho_{2}})\,
dy_{2}\leqslant\\\\
\displaystyle
\leqslant a\|\rho_{2}^{-1/2}\|_{C[0, d_{2}]}\|\partial_{y_{2}}\ln{\rho_{2}}\|_{L_{2}(0, d_{2})}\|\rho_{2}^{-1/2}(u_{1}-u_{2})\|_{L_{2}(0, d_{2})},
\end{array}
\end{equation}
\begin{equation}\label{step4_new1}\begin{array}{c}\displaystyle
\mu_{21}\int\limits^{d_{2}}_{0}(\partial_{y_{2}}(\rho_{2}\partial_{y_{2}}u_{1}))(\partial_{y_{2}}\ln{\rho_{2}})\,dy_{2}\leqslant\\ \\
\displaystyle
\leqslant |\mu_{21}|\,\|\rho_{2}^{-1/2}\|_{C[0, d_{2}]}\|\rho_{2}^{1/2}\partial_{y_{2}}(\rho_{2}\partial_{y_{2}}u_{1})\|_{L_{2}(0, d_{2})}\|\partial_{y_{2}}\ln{\rho_{2}}\|_{L_{2}(0, d_{2})}.
\end{array}
\end{equation}
It follows from \eqref{postzad_1}, \eqref{postzad_3} and \eqref{postzad_4} that for each $t\in[0,T]$ in at least one point $\varsigma_{2}(t)\in[0, d_{2}]$ the following relation holds
\begin{equation}\label{step4_8}
\rho_{2}(\varsigma_{2}(t), t)=d_{2}.
\end{equation}
Hence
\begin{equation}\label{step4_9}
\rho_{2}^{-1/2}(y_{2},t)=d_{2}^{-1/2}+\int\limits_{\varsigma_{2}}^{y_{2}}\partial_{s}\rho_{2}^{-1/2}(s,t)\, ds,
\end{equation}
and by the Cauchy inequality we have
\begin{equation}\label{step4_10}
\|\rho_{2}^{-1/2}\|_{C[0, d_{2}]}\leqslant d_{2}^{-1/2}+\frac{1}{2}\|\partial_{y_{2}}\ln{\rho_{2}}\|_{L_{2}(0, d_{2})}.
\end{equation}
Thus after integration of \eqref{step4_5} in $t$, using \eqref{postzad_10}, \eqref{step4_6}--\eqref{step4_new1} and \eqref{step4_10}, we derive
\begin{equation}\label{step4_11}\begin{array}{c}\displaystyle
\mu_{22}\|\partial_{y_{2}}\ln{\rho_{2}}\|^{2}_{L_{2}(0, d_{2})}+2K_{2}\gamma_{2}\int\limits_{0}^{t}\int\limits^{d_{2}}_{0}\rho_{2}^{\gamma_{2}}|\partial_{y_{2}}\ln{\rho_{2}}|^{2}\,
dy_{2}d\tau\leqslant\\\\
\displaystyle
\leqslant \mu_{22}\|\partial_{y_{2}}\ln{\widetilde{\widetilde{\rho}}_{02}}\|^{2}_{L_{2}(0, d_{2})} -2\int\limits^{d_{2}}_{0}u_{2}(\partial_{y_{2}}\ln{\rho_{2}})\,
dy_{2}+\\ \\
\displaystyle
+2\int\limits^{d_{2}}_{0}\widetilde{\widetilde{u}}_{02}(\partial_{y_{2}}\ln{\widetilde{\widetilde{\rho}}_{02}})\,
dy_{2}+2\int\limits_{0}^{t}\int\limits^{d_{2}}_{0}\rho_{2}|\partial_{y_{2}}u_{2}|^{2}\,
dy_{2}d\tau+\\ \\
\displaystyle
+|\mu_{21}|\!\int\limits_{0}^{t}\!\|\rho_{2}^{1/2}\partial_{y_{2}}(\rho_{2}\partial_{y_{2}}u_{1})\|_{L_{2}(0, d_{2})}\!\Big(2d_{2}^{-1/2}\!+\\\\
\displaystyle
+\|\partial_{y_{2}}\ln{\rho_{2}}\|_{L_{2}(0, d_{2})}\Big)\!\|\partial_{y_{2}}\ln{\rho_{2}}\|_{L_{2}(0, d_{2})} d\tau+\\ \\
\displaystyle
+a\int\limits_{0}^{t}\|\rho_{2}^{-1/2}(u_{1}-u_{2})\|_{L_{2}(0, d_{2})}\Big(2d_{2}^{-1/2}+\\\\
\displaystyle
+\|\partial_{y_{2}}\ln{\rho_{2}}\|_{L_{2}(0, d_{2})}\Big)\|\partial_{y_{2}}\ln{\rho_{2}}\|_{L_{2}(0, d_{2})}\,d\tau.
\end{array}
\end{equation}
Using the Cauchy inequality
$$
-2\int\limits^{d_{2}}_{0}u_{2}(\partial_{y_{2}}\ln{\rho_{2}})\,
dy_{2}\leqslant \frac{\mu_{22}}{2}\|\partial_{y_{2}}\ln{\rho_{2}}\|^{2}_{L_{2}(0, d_{2})}+\frac{2}{\mu_{22}}\|u_{2}\|^{2}_{L_{2}(0, d_{2})}
$$
and the estimates (see \eqref{step3_5})
\begin{equation}\label{step4_new2}
\|\partial_{xx}u_{1}\|^{2}_{L_{2}(Q_{T})}=\int\limits_{0}^{T}\int\limits_{0}^{d_{2}}\rho_{2}|\partial_{y_{2}}(\rho_{2}\partial_{y_{2}}u_{1})|^{2} \,dy_{2}d\tau\leqslant B_{12}^2
\end{equation}
and \eqref{step1_7}, from \eqref{step4_11} we deduce
\begin{equation}\label{step4_12}\begin{array}{c}\displaystyle
\|\partial_{y_{2}}\ln{\rho_{2}}\|^{2}_{L_{2}(0, d_{2})}\leqslant B_{14}+B_{15}\int\limits_{0}^{t}\Big(\|\rho_{2}^{-1/2}(u_{1}-u_{2})\|_{L_{2}(0, d_{2})}+\\ \\
\displaystyle
+\|\rho_{2}^{1/2}\partial_{y_{2}}(\rho_{2}\partial_{y_{2}}u_{1})\|_{L_{2}(0, d_{2})}\Big)\|\partial_{y_{2}}\ln{\rho_{2}}\|^{2}_{L_{2}(0, d_{2})}\,d\tau,
\end{array}
\end{equation}
where $B_{14}=B_{14}\left(M_{0}, B_{1}, B_{12}, \mu_{21}, \mu_{22}, d_{2}, a, T, \|\partial_{y_{2}}\ln{\widetilde{\widetilde{\rho}}_{02}}\|_{L_{2}(0, d_{2})},\right.$\linebreak $\left.\|\widetilde{\widetilde{u}}_{02}\|_{L_{2}(0, d_{2})}\right)$,
$B_{15}=B_{15}(\mu_{21}, \mu_{22}, a)$. Using the estimates \eqref{step1_7} and \eqref{step4_new2} again, by the Gronwall lemma from \eqref{step4_12} we get
\begin{equation}\label{step4_13}
\|\partial_{y_{2}}\ln{\rho_{2}}(t)\|_{L_{2}(0, d_{2})}\leqslant B_{16}\quad \forall\, t\in [0, T],
\end{equation}
where $B_{16}=B_{16}(B_{1}, B_{12}, B_{14}, B_{15}, a, T)$. Hence,
\begin{equation}\label{step4_14}\begin{array}{c}\displaystyle
|\ln\rho_{2}(y_{2},t)|=\left|\ln\rho_{2}(\varsigma_{2},t)+\int\limits_{\varsigma_{2}}^{y_{2}}\partial_{s}\ln\rho_{2}(s,t)\, ds\right|\leqslant \\ \\
\displaystyle
\leqslant |\ln{d_{2}}|+\sqrt{d_{2}}\|\partial_{y_{2}}\ln{\rho_{2}}\|_{L_{2}(0, d_{2})}\leqslant B_{17}(B_{16}, d_{2}),
\end{array}
\end{equation}
and consequently
\begin{equation}\label{step4_15}
0<B_{18}^{-1}\leqslant\rho_{2}(y_{2}, t)\leqslant B_{18}<\infty\quad\text{as}\quad(y_{2},t)\in[0, d_{2}]\times[0, T],
\end{equation}
where $B_{18}=e^{B_{17}}$. As a result,
\begin{equation}\label{step4_16}
0<B_{18}^{-1}\leqslant\rho_{2}(x, t)\leqslant B_{18}<\infty\quad\text{as}\quad(x,t)\in[0, 1]\times[0, T].
\end{equation}

\section{Estimates for the Derivatives of $\rho_{2}$ and $u_{2}$}
\label{sec7}

\noindent\indent It follows from \eqref{step4_13} and  \eqref{step4_16} that
\begin{equation}\label{step5_1}
\|\partial_{x}\rho_{2}(t)\|_{L_{2}(0, 1)}\leqslant B_{19}=B_{16}B_{18}^{\frac{3}{2}}\quad \forall\, t\in [0, T].
\end{equation}
We square the equation \eqref{postzad_2} with the number $ i = 2 $ and obtain
\begin{equation}\label{step5_2}\begin{array}{c}\displaystyle
\rho_{2}|\partial_{t}u_{2}|^{2}+\frac{\mu_{22}^{2}}{\rho_{2}}|\partial_{xx}u_{2}|^{2}-2\mu_{22}(\partial_{t}u_{2})(\partial_{xx}u_{2})=\\ \\
\displaystyle
=\frac{1}{\rho_{2}}\left|\mu_{21}\partial_{xx}u_{1}-K_{2}\partial_{x}\rho_{2}^{\gamma_{2}}-\rho_{2}u_{2}\partial_{x}u_{2}+a(u_{1}-u_{2})\right|^{2}.
\end{array}
\end{equation}
Let us consider the function
$$
\beta(t)=\mu_{22}\int\limits_{0}^{1}|\partial_{x}u_{2}|^{2}\, dx+\int\limits_{0}^{t}\int\limits_{0}^{1}\left(\rho_{2}|\partial_{t}u_{2}|^{2}+\frac{\mu_{22}^2}{\rho_{2}}|\partial_{xx}u_{2}|^{2}\right)\, dxd\tau.
$$
It follows from \eqref{postzad_4}, \eqref{step5_2} and the estimates \eqref{step4_16}, \eqref{step5_1} that
\begin{equation}\label{step5_3}
\begin{array}{c}
\displaystyle
\beta^{\prime}(t)\leqslant B_{20}+B_{21}\left(\|u_{1}-u_{2}\|^{2}_{L_{2}(0,1)}+\|\partial_{xx}u_{1}\|^{2}_{L_{2}(0,1)}\right)+\\\\
\displaystyle
+B_{22}\|u_{2}\|^{2}_{L_{\infty}(0,1)}\beta(t),
\end{array}
\end{equation}
where $B_{20}=B_{20}(B_{18}, B_{19}, K_{2}, \gamma_{2})$, $B_{21}=B_{21}(B_{18}, a, \mu_{21})$,\linebreak $B_{22}=B_{22}(B_{18}, \mu_{22})$, and in view of \eqref{step1_5}, \eqref{step2_17} and \eqref{step3_5} by the Gronwall lemma we get
\begin{equation}\label{step5_4}
\beta(t)\leqslant B_{23}\left(B_{1}, B_{12}, B_{20}, B_{21}, B_{22}, M_{0}, a, T, \mu_{22}, \|\partial_{x}u_{02}\|_{L_{2}(0, 1)}\right),
\end{equation}
and finally we arrive at the estimate
\begin{equation}\label{step5_5}
\begin{array}{c}
\|\partial_{x}u_{2}\|_{L_{\infty}\left(0, T; L_{2}(0, 1)\right)}+\|\partial_{xx}u_{2}\|_{L_{2}(Q_{T})}+\\\\
\displaystyle
+\|\partial_{t}u_{2}\|_{L_{2}(Q_{T})}\leqslant B_{24}(B_{18}, B_{23}, \mu_{22}).
\end{array}
\end{equation}
Now from the continuity equation \eqref{postzad_1} with the number $ i = 2 $ and the estimates \eqref{step4_16}, \eqref{step5_1}, \eqref{step5_5} it follows that
\begin{equation}\label{step5_6}
\|\partial_{t}\rho_{2}\|_{L_{\infty}(0, T; L_{2}(0,1))}\leqslant B_{25}(B_{18}, B_{19}, B_{24}).
\end{equation}
Thus, all estimates in the classes \eqref{postzad_5} for $ (\rho_{1}, \rho_{2}, u_{1}, u_{2}) $ are obtained for arbitrary $T\in (0, +\infty)$.

\section{Uniqueness of Solution}
\label{sec8}

\noindent\indent Let $\big(\rho^{(1)}_{1},\rho^{(1)}_{2}, u^{(1)}_{1}, u^{(1)}_{2}\big)$ and $\big(\rho^{(2)}_{1},\rho^{(2)}_{2}, u^{(2)}_{1}, u^{(2)}_{2}\big)$ be two strong (i.~e. of the class \eqref{postzad_5}) solutions to the problem \eqref{postzad_1}--\eqref{postzad_4}. Denote $\rho_{i}=\rho^{(1)}_{i}-\rho^{(2)}_{i}$,\linebreak $u_{i}=u^{(1)}_{i}-u^{(2)}_{i}$, $i=1,2$.

The equations \eqref{postzad_1} and initial data \eqref{postzad_3} yield the equalities
\begin{equation}\label{mamprok_ravlema52}
\partial_{t}\rho_{i}+\partial_{x}\left(\rho_i u_{i}^{(1)}\right)+\partial_{x}\left(\rho_{i}^{(2)} u_{i}\right)=0,\quad \rho_{i}|_{t=0}=0,\quad i=1,2.
\end{equation}
Multiplying \eqref{mamprok_ravlema52} by $2\rho_{i}$, and integrating in $x$ from $0$ to $1$, we obtain using \eqref{postzad_4} that $(i=1,2)$
\begin{equation}\label{mamprok_ravlema53}
\eta_{1i}^{\prime}(t)=-\int\limits^{1}_{0}\left(\rho_{i}^{2}\left(\partial_{x}u_{i}^{(1)}\right)+2\rho_{i}^{(2)}\rho_{i}(\partial_{x}u_{i})+
2\rho_{i}u_{i}\left(\partial_{x}\rho_{i}^{(2)}\right)\right)\, dx,
\end{equation}
where
\begin{equation}\label{mamprok_ravlemanew51}\eta_{1i}(t)=\int\limits^{1}_{0}\rho_{i}^{2}\, dx,\quad i=1, 2.\end{equation}
The terms in the right-hand side of \eqref{mamprok_ravlema53} are estimated as follows ($i=1,2$):
$$-\int\limits^{1}_{0}\rho_{i}^{2}\left(\partial_{x}u_{i}^{(1)}\right)\,dx\leqslant \|\partial_{x}u_{i}^{(1)}\|_{L_{\infty}(0, 1)}\eta_{1i}(t),$$
$$-2\int\limits^{1}_{0}\rho_{i}^{(2)}\rho_{i}(\partial_{x}u_{i})\,dx\leqslant\|\partial_{x}u_{i}\|^{2}_{L_{2}(0, 1)}+\|\rho_{i}^{(2)}\|_{L_{\infty}(Q_{T})}^{2}\eta_{1i}(t),$$
$$-2\int\limits^{1}_{0}\rho_{i}u_{i}\left(\partial_{x}\rho_{i}^{(2)}\right)\,dx\leqslant \|\partial_{x}\rho_{i}^{(2)}\|_{L_{\infty}\big(0, T; L_{2}(0, 1)\big)}^{2}\|u_{i}\|^{2}_{L_{\infty}(0, 1)}+\eta_{1i}(t)\leqslant$$
$$\leqslant \|\partial_{x}\rho_{i}^{(2)}\|_{L_{\infty}\big(0, T; L_{2}(0, 1)\big)}^{2}\|\partial_{x}u_{i}\|^{2}_{L_{2}(0, 1)} +\eta_{1i}(t).$$
Due to the inclusions
\begin{equation}\label{vkl050418}
\begin{array}{c}
\displaystyle
\rho_{i}^{(2)}\in L_{\infty}(Q_{T}),\quad \partial_{x}\rho_{i}^{(2)}\in L_{\infty}\big(0, T; L_{2}(0, 1)\big),\quad i=1,2,\\\\
\displaystyle
\partial_{x}u_{i}^{(1)}\in L_{2}\big(0, T; L_{\infty}(0, 1)\big),\quad i=1,2,
\end{array}
\end{equation}
we derive the estimates 
\begin{equation}\label{mamprok_ravlemanew52}\eta_{1i}^{\prime}(t)\leqslant B_{26}(t)\eta_{1i}(t)+B_{27}\|\partial_{x}u_{i}\|^{2}_{L_{2}(0, 1)},\quad i=1,2,\end{equation}
where $B_{26}\in L_{2}(0,T)$, while $B_{27}$ and $\|B_{26}\|_{L_{2}(0, T)}$ depend on $T$ and the norms in \eqref{vkl050418}. Hence, since $\eta_{1i}(0)=0$, $i=1,2$
(due to the initial data in ~\eqref{mamprok_ravlema52}), then applying the Gronwall inequality, we arrive at the inequalities
\begin{equation}\label{mamprok_ravlema54}
\eta_{1i}(t)\leqslant B_{28}\|\partial_{x}u_{i}\|^{2}_{L_{2}(Q_{t})},\quad i=1,2,
\end{equation}
where $B_{28}=B_{28}\left(B_{27}, T, \|B_{26}\|_{L_{2}(0, T)}\right)$.

Then, the equations \eqref{postzad_1}, \eqref{postzad_2} and the boundary conditions \eqref{postzad_4} imply the equality
\begin{equation}\label{mamprok_ravlema55}\begin{array}{c}\displaystyle
\frac{1}{2}\sum\limits_{i=1}^{2}\frac{d}{dt}\int\limits^{1}_{0}\rho_{i}^{(1)}u_{i}^{2}\,dx+
\sum\limits_{i, j=1}^2 \mu_{ij}\int\limits_{0}^{1}(\partial_{x}u_{i})(\partial_{x}u_{j})\, dx+
\\ \\
\displaystyle +a\int\limits_{0}^{1}|u_{1}-u_{2}|^{2}\, dx=\sum\limits_{i=1}^{2}K_{i}\int\limits_{0}^{1}\left(\left(\rho^{(1)}_{i}\right)^{\gamma_{i}}-\left(\rho^{(2)}_{i}\right)^{\gamma_{i}}\right)(\partial_{x}u_{i})\, dx-\\ \\
\displaystyle -\sum\limits_{i=1}^{2}\int\limits_{0}^{1}\rho_{i}u_{i}\left(\partial_{t}u_{i}^{(2)}\right)\,dx-
\sum\limits_{i=1}^{2}\int\limits_{0}^{1}\rho_{i}^{(1)}u_{i}^{2}\left(\partial_{x}u_{i}^{(2)}\right)\,dx-\\\\
\displaystyle
-\sum\limits_{i=1}^{2}\int\limits_{0}^{1}\rho_{i}u_{i}u_{i}^{(2)}\left(\partial_{x}u_{i}^{(2)}\right)\,dx,
\end{array}\end{equation}
from which it follows that
\begin{equation}\label{mamprok_ravlemanew53}\begin{array}{c}\displaystyle \eta_{2}^{\prime}(t)+\eta_{3}^{\prime}(t)\leqslant \sum\limits_{i=1}^{2}
K_{i}\int\limits_{0}^{1}\left(\left(\rho^{(1)}_{i}\right)^{\gamma_{i}}-\left(\rho^{(2)}_{i}\right)^{\gamma_{i}}\right)(\partial_{x}u_{i})\, dx
-\\\\
\displaystyle
-\sum\limits_{i=1}^{2}\int\limits_{0}^{1}\rho_{i}u_{i}\left(\partial_{t}u_{i}^{(2)}\right)\,dx-
\sum\limits_{i=1}^{2}\int\limits_{0}^{1}\rho_{i}^{(1)}u_{i}^{2}\left(\partial_{x}u_{i}^{(2)}\right)\,dx-\\ \\
\displaystyle
-\sum\limits_{i=1}^{2}\int\limits_{0}^{1}\rho_{i}u_{i}u_{i}^{(2)}\left(\partial_{x}u_{i}^{(2)}\right)\,dx,
\end{array}\end{equation}
where
\begin{equation}\label{mamprok_ravlemanew54}\eta_{2}(t)=\frac{1}{2}\sum\limits_{i=1}^{2}\int\limits_{0}^{1}\rho_{i}^{(1)}u_{i}^{2}\,dx,\quad \eta_{3}(t)=M_{0}\sum\limits_{i=1}^{2}\int\limits_{0}^{t}\int\limits_{0}^{1}|\partial_{x}u_{i}|^{2}\,dx d\tau.
\end{equation}
We estimate the terms in the right-hand side of the relation \eqref{mamprok_ravlemanew53} as follows:
$$\sum\limits_{i=1}^{2}K_{i}\int\limits_{0}^{1}\left(\left(\rho^{(1)}_{i}\right)^{\gamma_{i}}-\left(\rho^{(2)}_{i}\right)^{\gamma_{i}}\right)(\partial_{x}u_{i})\, dx \leqslant \frac{M_{0}}{4}\sum\limits_{i=1}^{2}\|\partial_{x} u_{i}\|_{L_{2}(0, 1)}^{2}+$$
$$+B_{29}\left(M_{0}, \{K_{i}\}, \{\gamma_{i}\}, \left\{\|\rho_{i}^{(1, 2)}\|_{L_{\infty}(Q_{T})}\right\}\right)\sum\limits_{i=1}^{2}\eta_{1i}(t),$$
$$-\sum\limits_{i=1}^{2}\int\limits_{0}^{1}\rho_{i}u_{i}\left(\partial_{t}u_{i}^{(2)}\right)\,dx\leqslant$$
$$\leqslant\frac{M_{0}}{4}\sum\limits_{i=1}^{2}
\|\partial_{x} u_{i}\|_{L_{2}(0, 1)}^{2}+B_{30}\left(M_{0}, \left\{\|\partial_{t}u_{i}^{(2)}\|_{L_{2}(0, 1)}\right\}\right)\sum\limits_{i=1}^{2}\eta_{1i}(t),$$
$$-\sum\limits_{i=1}^{2}\int\limits_{0}^{1}\rho_{i}^{(1)}u_{i}^{2}\left(\partial_{x}u_{i}^{(2)}\right)\,dx\leqslant B_{31}
\left(\left\{\|\partial_{x}u_{i}^{(2)}\|_{L_{\infty}(0, 1)}\right\}\right)\eta_{2}(t),$$
$$-\sum\limits_{i=1}^{2}\int\limits_{0}^{1}\rho_{i}u_{i}u_{i}^{(2)}\left(\partial_{x}u_{i}^{(2)}\right)\,dx\leqslant B_{32}\left(\left\{\|u_{i}^{(2)}\|_{L_{\infty}(Q_{T})}\right\}\right)\sum\limits_{i=1}^{2}\eta_{1i}(t)+$$
$$+B_{33}\left(\left\{\|\partial_{x}u_{i}^{(2)}\|_{L_{\infty}(0, 1)}\right\}, \left\{\|1/\rho_{i}^{(1)}\|_{L_{\infty}(Q_{T})}\right\}\right)\eta_{2}(t),$$
where $B_{30}, B_{31}, B_{33}\in L_{1}(0, T)$. Thus, from \eqref{mamprok_ravlemanew53}, taking into account the already proved relation (see \eqref{mamprok_ravlema54})
$$\sum\limits_{i=1}^{2}\eta_{1i}(t)\leqslant B_{34}(M_{0}, B_{28})\eta_{3}(t),$$
we derive that
\begin{equation}\label{mamprok_ravlema513}\eta_{2}^{\prime}(t)+\frac{1}{2}\eta_{3}^{\prime}(t)\leqslant B_{35}(t)\left(\eta_{2}(t)+\frac{1}{2}\eta_{3}(t)\right),\end{equation}
where $B_{35}=B_{31}+B_{33}+2B_{34}(B_{29}+B_{30}+B_{32})\in L_{1}(0, T)$, and finally, due to $\eta_{2}(0)=\eta_{3}(0)=0$, we arrive at the identities
\begin{equation}\label{mamprok_ravlema514}\eta_{11}\equiv \eta_{12}\equiv \eta_{2}\equiv \eta_{3}\equiv 0,\end{equation}
which complete the proof of Theorem 2.

\end{document}